\documentclass[12pt,a4paper]{article}
\usepackage{t1enc}
\usepackage[latin2]{inputenc} 
\usepackage[dvips,final]{graphicx} 
\usepackage{tabularx} 
\usepackage{amsfonts}
\usepackage{amsthm}
\usepackage{amssymb}
\usepackage{epsf}

\def\du{\unskip\smash{\lower 1.4ex \hbox{\char34}}\kern-.2ex}
\def\hu{\kern-.2ex\hbox{\char92}}

\newcommand{\bdis}{\begin{displaymath}}
\newcommand{\edis}{\end{displaymath}}
\newcommand{\be}{\begin{equation}}
\newcommand{\ee}{\end{equation}}
\newcommand{\ep}{\epsilon}

\def\Diff{\mbox{\it Diff }}
\def\diff{\mbox{\it diff }}
\def\Der{\mbox{\it Der }}
\def\deg{\mbox{\it deg }}
\def\Pol{\mbox{ Pol }}
\def\M{\mfrak{M}}
\def\R{\mbb{R}^{0|n}}
\def\Rm{\mbb{R}^{0|2n}}

\def\Rn{\mbb{R}^{n}}
\def\Rz{\mbb{R}^{0|1}}
\def\X{\mfrak{X}}
\def\C{C^\infty}

\def\ua{\uparrow}

\def\a{\alpha}
\def\b{\beta}
\def\g{\gamma}
\def\m{\mu}
\def\n{\nu}
\def\l{\lambda}
\def\T_n{T_{{}\nabla}}
\def\U_n{R_{{}\nabla}}

\newcommand{\ra}{\rightarrow}

\newcommand{\mcal}{\mathcal}
\newcommand{\mfrak}{\mathfrak}
\newcommand{\pd}{\partial}
\newcommand{\mbb}{\mathbb}

\theoremstyle{citing}

\begin{document}
\title{{\textbf{Grassmann Electrodynamics and General Relativity}}}
\author{Denis Kochan \\
{\it{Department of Theoretical Physics}}, \\ {\it{Faculty of
Mathematics Physics and Informatics}},
\\ {\it{ Comenius University}},\\ {\it{ Mlynsk\' a Dolina F2, 842 48 Bratislava, Slovakia}} \\
 \tt{e-mail: kochan@sophia.dtp.fmph.uniba.sk} }
\date{}
\maketitle
{\abstract{\noindent The aim of this paper is to present a short introduction to
supergeometry on pure odd supermanifolds. (Pseudo)differential forms, Cartan calculus
(DeRham differential, Lie derivative, "inner" product), metric, inner product,
Killing's vector fields, Hodge star operator, integral forms, co-differential and
connection on odd Riemannian supermanifolds are introduced.
The electrodynamics and Einstein relativity with anti-commuting variables only are
formulated modifying the geometry beyond classical
(even, bosonic) theories appropriately.
Extension of these ideas to general supermanifolds is straightforward.
}}

\section{Introduction}

Supergeometry is an interesting and fruitful branch of mathematics
with a variety of powerful applications in modern theoretical physics,
in particular in SUSY, supergravity and superstrings. From a purely mathematical
point of view, supergeometry is natural extension of the ordinary
differential geometry by Grassmann variables. Such anti-commuting extensions
represent an essential and inspiring feature of all supermathematics.\\
The first paper about supermathematics was the work of
Martin \cite{Martin}, in which the classical limit of a system with fermionic
degrees of freedom was discussed. This theory, later called pseudoclassical mechanics
was independently developed in the middle of 70thies
by Berezin and Marinov \cite{Berezin-Marinov1}-\cite{Berezin-Marinov2}, Casalbuoni
\cite{Casalbuoni1}-\cite{Casalbuoni2} and others. Since the Grassmann variables
became an invaluable tool in the description of fermions, and because their natural combination
with even (bosonic) degrees of freedom led at the beginning of 70thies to the discovery of
supersymmetry, it was a necessary to build a rigorous mathematical theory, which would be
able to describe both (even and odd) degrees of freedom. Systematical investigation in this
direction was initiated at the beginning of 60thies by Berezin \cite{Berezin1}-\cite{Berezin2},
but the main goals of the supermathematics were established during 70thies largely by the
Russian mathematical school led by Berezin. More details about supergeometry
(and also about its chronology) can be found in the review article \cite{Leites},
and in the famous Berezin book \cite{Berezin} which could be indeed regarded as the Bible of
supermathematics (see also references therein).\\
The aim of this paper is to present a very short introduction to supergeometry over
pure odd supermanifolds (sections 2-4). Using supergeometrical methods, we
explain the electromagnetism described only by anti-commuting coordinates (section 5);
such extravagant theory is called with a grain of salt \emph{Grassmann electrodynamics}.
After the definition of linear connection on pure odd Riemannian supermanifold
we shall be able to reveal an odd analogy of the resulting theory with the Einstein
theory of relativity (section 6).
All this "odd business" (in both meanings of the word "odd") is based on classical geometrical
analogy, similarly as pseudoclassical mechanics developed by Martin, Berezin and Casalbuoni.

\section{(Pseudo)differential forms}

The (pseudo)differential forms on an arbitrary smooth (real) $m|n$-dimensional supermanifold
$\M$ was in general investigated in the framework which is applied below in \cite{Kochan}.\\
We shall study the basic properties of the (pseudo)dif\-fe\-ren\-tial
forms on a pure odd ($0|n$-dimensional) real supermanifold, i.e. on $\R$.
The odd $n$-dimensional Cartesian space $\R$ is covered by $n$ global Grassmann coordinates
$(\xi^1,\dots ,\xi^n)$ and the superalgebra of functions $\mcal{F}(\R)$
coincides with exterior ($\mbb{Z}_{2}$-graded) algebra
$\bigwedge\Rn=[\bigwedge\Rn]_{[0]}\oplus[\bigwedge\Rn]_{[1]}$.\\
The tangent bundle of $\R$ is a supermanifold $T\R=\R\times\R$ with a set of global anti-commuting
coordinates $(\xi^1,\dots ,\xi^n,\sigma^1,\dots ,\sigma^n)$ transforming under the
transformation of coordinates $\xi^\a\mapsto \Xi^\a(\xi)$ on the base $\R$ as
\be\label{tangent bundle}
(\xi^\a,\,\sigma^\a) \mapsto \biggl(\Xi^\a(\xi),\,\Sigma^\a(\xi,\sigma)=\sigma^\b\frac{\pd\, \Xi^\a}{\pd\, \xi^\b}\biggr)\ .
\ee
The odd functions $\Xi^\a(\xi)$ in (\ref{tangent bundle}) guarantee that the parity
of coordinates on $T\R$ is preserved.
Throughout the paper, we use left derivatives with respect to Grassmann variables
and Einstein summation convention. The parity of any object $O$ (with respect to anti-commuting variables)
is denoted by $\tilde{O}$, and to distinguish Grassmann and ordinary variables we use Greek letters
for the former and Latin letters for the latter.\\
The odd tangent bundle $\Pi T\R=\R\times\mbb{R}^n$ is a supermanifold which is obtained from $T\R$ by changing
the parity of the fiber variables $\sigma^\a$. The coordinate transformation on supermanifold $\R$
induces the corresponding transformation on odd tangent bundle
\be\label{odd tangent bundle}
(\xi^\a,\,y^\a) \mapsto \biggl(\Xi^\a(\xi),\,Y^\a(\xi,y)=y^\b\frac{\pd\, \Xi^\a}{\pd\, \xi^\b}\biggr)\ .
\ee
It is well known (for more details see \cite{Manin}, but we hope that it will become clear from our next explanation)
that the superalgebra of {\textbf{\emph{differential forms}} on $\R$ can be identified with
$\mbb{Z}_2$-graded algebra $\Pol(\Pi T\R)=\Pol(\mbb{R}^n)\otimes\bigwedge\Rn$
of all polynomials with real coefficients over supermanifold $\Pi T\R$.
There is one to one correspondence between the differentials of Grassmann variables ${\rm d}\xi^\a$
(even quantities) and the even variables $y^\a$ covering fibers in $\Pi T\R$. Their natural generalization
leads to the definition of {\textbf{\emph{pseudodifferential forms}} over $\R$, namely, the superalgebra of
pseudodifferential forms is defined as the $\mbb{Z}_2$-graded algebra $\C(\Pi T\R)=\C(\mbb{R}^n)\otimes\bigwedge\Rn$.\\
The standard differential operations on forms, DeRham differential, Lie derivative and
inner product, are identified with special vector fields on $\Pi T\R$. To obtain their exact forms we use the
fruitful idea of Maxim Kontsevich, who pointed out (see \cite{Kontsevich}) that
the odd tangent bundle of arbitrary $m|n$-dimensional supermanifold $\M$
is canonically isomorphic to the supermanifold of all supermaps $\Rz\rightarrow\M$.
In our case
\be
\Pi T\R\equiv\{\mbox{\ supermaps:\ }\Rz\rightarrow\R\}\ .
\ee
An arbitrary supermap $\Phi\in\Pi T\R$ is expressed in coordinates
(by using Taylor expansion in $\theta$) as
\bdis
\Phi:\ \theta\mapsto\xi^\a(\Phi(\theta))=\xi^\a+\theta y^\a\ .
\edis
It is clear that such $\Phi$ is characterized by $n$ odd and $n$ even
coordinates, which transform in accordance with (\ref{odd tangent bundle}).\\
The supergroup $\Diff(\Rz)=\{\mbox{\,diffeomorphisms:\,} \Rz\ra\Rz\,;\theta\mapsto\theta^\prime=\theta a+\b\}$
defines via its natural right action
\bdis
\Pi T\R\times\Diff(\Rz)\ra\Pi T\R\ , \ \ \ \ \ \ \ \ \ \ \  (\Phi,g)\mapsto \Phi\circ g\ ,
\edis
the left invariant (fundamental) vector fields $E,\,Q$ on $\Pi T\R$. Their expression
in coordinates is very simple, namely
\begin{eqnarray}
E & = & y^\a\pd_{y^\a}\ \ \ \ \mbox{\emph{\textbf{Euler field}}}\ \ (\tilde{E}=0)\ ,\label{Euler}\\
Q & = & y^\a\pd_{\xi^\a} \ \ \ \ \mbox{\emph{\textbf{DeRham differential}}}\ \ (\tilde{Q}=1)\ .\label{deRham}
\end{eqnarray}
The Euler vector field "measures" the degree of homogenity of (pseudo)differential forms under the
supergroup action, therefore the superalgebra $\C(\Pi T\R)$ has also a natural
$\mbb{Z}$-graded structure ($f\in\left[\C(\Pi T\R)\right]^{(k)}$ $\Leftrightarrow$ $Ef=kf=:\deg(f)f$).
A direct calculation gives the (super)commutation relations in the Lie superalgebra $\diff(\Rz)$:
\be
[E,\, E]=0\ , \ \ \ \ \ \ \ [E,\, Q]=Q\ , \ \ \ \ \ \ \ [Q,\, Q]=2Q^2=0\ .
\ee
Similarly, the supergroup $\Diff(\R)=\{\mbox{\,diffeomorphisms of\ }\R\}$
acts on the odd tangent bundle $\Pi T\R$,
\bdis
\Diff(\R)\times\Pi T\R\ra\Pi T\R\ , \ \ \ \ \ \ \ \ \ \ \  (g,\Phi)\mapsto g\circ \Phi\ .
\edis
Therefore to any element $V=V(\xi^\a)\pd_{\xi^\a}=V^\a(\xi)\pd_{\xi^\a}$ of the
corresponding Lie superalgebra $\diff(\R)=\X(\R)=\Der(\mcal{F}(\R))$
we can assign unique vector field $V^\uparrow$ on $\Pi T\R$.
A straightforward coordinate computation\footnote{In the case of an odd vector field
($\tilde{V}=1$) it is necessary to consider instead of flow the {\emph{superflow}}
(homomorphism of the supergroups $\mbb{R}^{1|1}$
and $\Diff(\R)$), whose infinitesimal ($\Delta t,\Delta\ep$) action
in the coordinates is:
$
\xi^\a\mapsto\xi^\a+\Delta\ep V(\xi^\a)+\frac{\Delta t}{2}[V,\, V](\xi^\a)
$, $\Delta\ep$ is odd variable.}
gives
\be\label{Lie derivation}
V^\ua=V(\xi^\a)\pd_{\xi^\a}+(-1)^{\tilde{V}}Q(V(\xi^\a))\pd_{y^\a}\ \ \Rightarrow\ \ \tilde{V^\ua}=\tilde{V}\ .
\ee
Apart from this natural lifting construction, it is also possible to associate to any $V\in\X(\R)$
certain vector field $V_\ua$ on $\Pi T\R$ such that $\tilde{V_\ua}=\tilde{V}+1$ and
\be\label{Cartan magic formula}
{}[V_\ua,\, Q]=V^\ua\ .
\ee
Obviously, the coordinate expression for $V_\ua$ is
\be\label{inner product}
V_\ua=V(\xi^\a)\pd_{y^\a}\ .
\ee
For any vector fields $V,\,W\in\X(\R)$ it is easy to confirm the validity of supercommutations relations
\be
\begin{array}{lcccl}
   {} [E,\, V^\ua]\ =0\ , & & & & [E,\, V_\ua]=-V_\ua\ , \\
   {} [V^\ua,\, Q]\ =0\ , & & & & [V^\ua,\, W^\ua]=[V,\, W]^\ua \ ,\\
  {}[V_\ua,\, W_\ua]=0\ , & & & & [V^\ua,\, W_\ua]\,=[V,\, W]_\ua\ .
\end{array}
\ee
The vector field $V^\ua$ corresponds to the {\textbf{\emph{Lie derivative}}} $\mcal{L}_V$ (with respect to $V$)
acting on forms, whereas $V_\ua$ represents the {\textbf{\emph{inner product}}} $i_V$ (with $V$).
Equation (\ref{Cartan magic formula}) is the famous {\emph{Cartan formula}}.\\
An arbitrary (pseudo)differential form is a polynomial (function) on the supermanifold $\Pi T\R$
and therefore it can be expressed in any local coordinates as
\be
f=f(\xi,y)=\sum\limits_{\beta=0}^{n}\sum\limits_{\a_1,\dots ,\a_\beta}f_{\a_1,\dots ,\a_\beta}(y)\wedge\xi^{\a_1}\wedge\dots\wedge\xi^{\a_\beta}\ ,
\ee
with ordinary real polynomials (functions) $f_{\a_1,\dots ,\a_\b}(y)$, which are
skew-sym\-metric in the indices ${\a_1,\dots ,\a_\b}$.\\
The {\emph{\textbf{integral}}} of the pseudodifferential form $f$ over $\R$ is defined as Berezin integral
(for more details see \cite{Berezin},\cite{Manin}) of a function $f$ on $\Pi T\R$:
\be\label{integral}
\mbox{\textbf{I}}[f]:=\int\limits_{\bigwedge\Rn}\overline{{\rm d}\xi}\int\limits_{\Rn}\overline{{\rm d}y}\ f(\xi,y)\ .
\ee
It is clear that such integral is not well defined for all elements of the superalgebras $\Pol(\Pi T\R)$
and $\C(\Pi T\R)$, because the manifold $\Rn$ (the typical fiber in $\Pi T\R$) is not compact.
The berezinian of the transformation (\ref{odd tangent bundle}) is equal to unity,
therefore the integral (\ref{integral}) is coordinate independent; moreover,
the odd \emph{per-partes} integration formula (the odd \emph{Stokes theorem}) holds,
\be\label{per-partes & Stokes theorem}
\mbox{\textbf{I}}[(Q f)\wedge g]=(-1)^{\tilde{f}+1}\,\mbox{\textbf{I}}[f \wedge(Q g)]\ \ \Leftrightarrow\ \ \mbox{\textbf{I}}[Q(f\wedge g)]=0\ .
\ee
An arbitrary supermap $\Phi:\R\rightarrow\mbb{R}^{0|m}$ defines (in accordance with (\ref{odd tangent bundle}))
the supermap $\Phi^{_\ua}:\Pi T\R\rightarrow\Pi T\mbb{R}^{0|m}$. The corresponding superalgebra homomorphism
${\Phi^{_\ua}}^*:\C(\Pi T\mbb{R}^{0|m})\rightarrow\C(\Pi T\R)$ is the \emph{pull-back} of the
supermap $\Phi$ over (pseudo)differential forms. For the vector field $V\in\X(\R)$, which generates the
(super)flow (diffeomorphism) on the supermanifold $\R$, the general formula for the \emph{pull-back} reads:
\be\label{general formula for superflow}
f\mapsto f_{(t,\,\ep)}\equiv[\Phi_{(t,\,\ep)}(V^\ua)]^*\,f=\biggl\{\begin{array}{ll}\exp{\{t\, V^\ua\}}\,f & \mbox{\ for\ } \tilde{V}=0 \ ,\\
\exp{\{\ep\, V^\ua+\frac{t}{2}\,[V^\ua,\,V^\ua]\}}\,f & \mbox{\ for\ }\tilde{V}=1\ .
\end{array}\biggr.
\ee

\section{Metric and Killing's vector fields}

The metric may be introduced on an arbitrary $m|n$-dimensional smooth supermanifold $\M$ (in particular,
on an ordinary manifold $M$) as an even regular (non-degenerate) quadratic function in fiber variables
on a tangent bundle $T\M$. In our case $\M=\R$ and the {\emph{\textbf{metric}}} has the form
\be\label{metric}
g=g(\xi,\sigma)=g_{\a\b}(\xi)\wedge\sigma^\a\wedge\sigma^\b\ ,
\ee
where functions $g_{\a\b}(\xi)=-g_{\b\a}(\xi)$ are even elements of $\mcal{F}(\R)$ (roughly speaking,
components of metric tensor in coordinates $(\xi^1,\dots ,\xi^n)$). The non-degeneracy
condition reads
\be\label{non-degeneracy condition}
|g|:=\det (g_{\a\b})\neq 0\ .
\ee
Let us emphasize that non-degeneracy of $g$ implies that the even skew-symmetric
matrix $g_{\a\beta}$ is invertible. Consequently, $0|n$-dimensional supermanifold $\R$
can be Riemannian only if $n$ is even (this fact is strongly reminiscent of the situation
in symplectic geometry), therefore our next analysis will be performed only for pure odd,
even-dimensional supermanifolds $\Rm$.\\
An arbitrary vector field $V=V(\xi^\a)\pd_{\xi^\a}\in\X(\Rm)$ could
be vertically lifted from the base supermanifold $\Rm$ to the tangent bundle $T\Rm$: the
coordinate expression for vertically lifted vector field over the tangent bundle is very simple
\be\label{vertical lift}
V_{\ua_{ver}}=V(\xi^\a)\pd_{\sigma^\a}\ \ \Rightarrow\ \ \tilde{V}_{\ua_{ver}}=\tilde{V}\ .
\ee
The metric $g$ on the supermanifold $\Rm$ allows us to define the {\emph{\textbf{inner product}}}
of vector fields on $\Rm$ as follows
\be\label{inner product on vector field}
(V\,,\,W)_g:=(-1)^{\tilde{W}+1}\,V_{\ua_{ver}}\biggl[W_{\ua_{ver}} \biggl(\frac{g}{2}\biggr)\biggr]=V(\xi^\a)\wedge g_{\a\b}(\xi)\wedge W(\xi^\b)\ .
\ee
It is clear that for all vector fields $V,\,W,\,U\in\X(\Rm)$ homogeneous with
respect to parity and arbitrary function $f\in\mcal{F}(\Rm)$ the following relations are valid
{\small{
\be\hspace{-0.4cm}
\begin{array}{l}
\begin{array}{rcl}
\hspace{1,12cm}\widetilde{(V\,,\,W)_g} & = & \tilde{V}+\tilde{W}\ ,\\
 (V\,,\,W)_g & = & -(-1)^{(\tilde{V}+1) (\tilde{W}+1)} (W\,,\,V)_g\ \ \mbox{odd graded skew-symmetry}\ ,
\end{array}\\
\left.\begin{array}{rcl}
(V+fU\,,\,W)_g & = & (V\,,\,W)_g+f(U\,,\,W)_g \\
(V\,,\,W+fU)_g & = & (V\,,\,W)_g+(-1)^{\tilde{V}\tilde{f}}f(V\,,\,U)_g
\end{array}\right\}\ \mbox{graded $f$-linearity}\ .
\end{array}
\ee}}%
In the ordinary differential geometry, it is well known, that apart from the vertical lifting procedure there also
exists a canonical horizontal lift of vector field from the base $\Rm$ to the tangent bundle $T\Rm$.
An arbitrary even [odd] vector field $V=V(\xi^\a)\pd_{\xi^\a}\in\X(\Rm)$ induces an infinitesimal flow [superflow]
on the base supermanifold $\Rm$. Such infinitesimal diffeomorphism of $\Rm$, in accordance with (\ref{tangent bundle}),
generates the (super)flow on the tangent bundle $T\Rm$. Its generator is the vector field
\be\label{horizontal lift}
V^{\ua_{hor}}=V(\xi^\a)\,\pd_{\xi^\a}+\Biggl[\sigma^\b\pd_{\xi^\b}\biggl(V(\xi^\a)\biggr)\Biggr]\pd_{\sigma^\a}\ \ \Rightarrow\ \ \tilde{V}^{\ua_{hor}}=\tilde{V} \ .
\ee
The origin of the horizontal lifted vector field $V^{\ua_{hor}}\in\X(T\Rm)$ is the same as the origin
of the vector field $V^{\ua}\in\X(\Pi T\Rm)$, which acts as \emph{Lie derivative} on the algebra of
(pseudo)differential forms. \\
The {\emph{\textbf{conformal Killing's vector fields}}} on the Riemannian
supermanifold $\Rm$  are solutions of the system of $n(2n-1)$ algebraic equations
{\small{
\be\label{Killing vector fields}
\chi\wedge g=V^{\ua_{hor}}(\,g\,)\Leftrightarrow \chi\wedge g_{\a\b}=V^\mu (g_{\a\b})_{,\,\mu}+(-1)^{\tilde{V}}\Biggl[g_{\a\mu}(V^\mu)_{,\,\b}-g_{\b\mu}(V^\mu)_{,\,\a}\Biggl]\ ,
\ee}}%
where $\chi\in\mcal{F}(\Rm)$ is even conformal scaling function and $(f)_{,\,\mu}=\pd_{\xi^\mu}[f(\xi)]$.
It is evident that the linear combination of two conformal Killing's vector fields is again the conformal
Killing's vector field, and because
\bdis
{}[V\,,\,W]^{\ua_{hor}}=[V^{\ua_{hor}}\,,\,W^{\ua_{hor}}]\ ,
\edis
the supercommutator of two conformal Killing's vector fields is a generator of conformal transformation of the
supermanifold $\Rm$, too.\\
It is possible to show that the Lie superalgebra of pure Killing vector fields ($\chi=0$) over $\Rm$ is at a most
$n(2n+1)|2n$-dimensional subsuperalgebra of $n2^{(2n)}|n2^{(2n)}$-dimensional $\mbb{Z}_2$-graded algebra $\X(\Rm)$.
The proof is analogical as in the ordinary differential geometry, but we do not prove this statement here, because
it is not necessary for our next construction and, moreover, it requires the definition of a new supergeometrical
notion, namely, the exponential (super)mapping.\\

\noindent Let us note that in analogy with ordinary differential geometry it is possible in supergeometry to define
objects, similarly as it was done with metric, which correspond (from ordinary geometrical
point of view) to covariant [contravariant] symmetric and anti-symmetric tensors:\\
\emph{Covariant} [\emph{contravariant}] \emph{symmetric tensor field} of rank $k$ over an arbitrary smooth
supermanifold $\M$ (in particular, on an ordinary manifold $M$) is defined as the polynomial function
of degree $k$ in fiber variables on the tangent [cotangent] bundle $T\M$ [$T^*\M$].\\
\emph{Anti-symmetric covariant} [\emph{contravariant}] \emph{tensors} (differential forms [multivector fields])
are analogically encoded in the polynomials in fiber variables on the odd tangent [odd cotangent]
bundle $\Pi T\M$ [$\Pi T^*\M$]. More detailed (but not exhaustive) description of the tensorial supercalculus
on smooth supermanifolds may be found, e.g. in \cite{Berezin},\cite{Manin},\cite{DeWitt}.

\section{Hodge $*_{g,\,o}\,$, Integral forms and Co-differential}

As in the case of ordinary differential forms, the metric is an essential ingredient in the definition of
Hodge $*_{g,\,o}$. Because now we are familiar with all its relevant ingredients, we are able
to define this operator.\\
For the function $f=f(\xi,y)\in[\C(\Pi T\Rm)]^{(k)}$ ((pseudo)differential forms over $\Rm$)
the {\textbf{\emph{Hodge $*_{g,\,o}$ operator}}} is formally defined by its Fourier transform in the fibre variables
$y^\a$, namely
\be\label{Hodge *}
(*_{g,\,o} f)(\xi,y):=\frac{(\imath)^{\hat{f}}}{(2\pi)^{n}}\int\limits_{\mbb{R}^{2n}}\overline{{\rm d}z}(o\, \sqrt{|g|})\wedge f(\xi,z)\wedge\exp{\{-\imath z^\a g_{\a\b} y^\b\}}\ ,
\ee
where the symbol $\hat{f}$ denotes the parity with respect to even variables, i.e $\hat{f}=k|\mbox{\,mod\,} 2$
and the orientation $o=\pm1$ (because the square root of $|g|$ is uniquely defined up to sign). In what follows,
to simplify notation, we will put $o=1$ and subscript $o$ will be omitted.
It is clear that the Hodge $*_g$ operator is defined on elements from $\Pol(\Pi T\Rm)$ only in sense of
distributions. Such generalized functions with one point support on the supermanifold $\Pi T\Rm$ are
called {\textbf{\emph{integral forms}}} over the odd tangent bundle and we denote them as  $*(\Pol(\Pi T\Rm))$.
A straightforward calculation shows that the definition of the Hodge star operator does not dependent on the
choice of coordinates. The definition (\ref{Hodge *}) is strictly correct only for functions from
$\C(\Pi T\Rm)$ that are behaving well in the variables $y^\a$ at infinity (e.g. functions with compact support).
The basic properties of the Hodge star operator can be obtained from the definition (\ref{Hodge *}):
\bdis
\widehat{*_g f}=\hat{f},\ \ \  \widetilde{*_g f}=\tilde{f},\ \ \  \widehat{Q(*_g f)}=\hat{f}+1,\ \ \  \widetilde{Q(*_g f)}=\tilde{f}+1,\ \ \  *_g(*_g f)=(-1)^{\hat{f}} f.
\edis
Similarly, like in the standard differential geometry, the metric $g$ and
the DeRham differential $Q$ define a new operator $\delta_g$ acting on (pseudo)differential
forms. This operator, called {\textbf{\emph{co-differential}}} is defined by equation
\be
Q(h)\wedge *_g f - Q(h\wedge *_g f)=:(-1)^{\tilde{h}}h\wedge *_g (\delta_g f)\ ,
\ee
where functions $h\,,f\in\C(\Pi T\Rm)$ are homogenous elements with respect
to Grassmann and fiber variables respectively. From this definition it follows
\be\label{co-differential}
\delta_g(f)=(-1)^{\hat{f}}*_g Q *_g f=\biggl[\frac{\pd\ }{\pd\,\xi^\a}-\frac{\pd\ln{\sqrt{|g|}}}{\pd\,\xi^\a}-y^\mu\,\frac{\pd g_{\mu\nu}}{\pd\,\xi^\a}\,g^{\nu\lambda}\frac{\pd\ }{\pd\,y^\lambda}\biggr]\,g^{\a\b}\frac{\pd f}{\pd\,y^\b}
\ee
where $g^{\a\mu}(\xi)$ is inverse to the matrix $g_{\mu\b}(\xi)$, i.e. $g^{\a\mu}(\xi)\wedge g_{\mu\b}(\xi)=\delta^\a_\b$.
All basic properties of the co-differential on Riemannian supermanifold $\Rm$ are given by equation
(\ref{co-differential}), and may be deduced from the properties of Hodge star operator and
DeRham differential. It is immediately evident from (\ref{Hodge *}) and (\ref{co-differential}) that
for arbitrary (pseudo)differential form $f$ and any supermap $\Phi:\Rm\rightarrow\mbb{R}^{0|2m}$
it holds
\be
{\Phi^{_\ua}}^*[*_g f]=*_{(\Phi^*g)}({\Phi^{_\ua}}^* f) \ \ \ \ \ \mbox{and}\ \ \ \ \ {\Phi^{_\ua}}^*[\delta_g (f)]=\delta_{(\Phi^*g)}({\Phi^{_\ua}}^* f)\ ,
\ee
where $\Phi^*g$ denotes the \emph{pull-back} of the metric $g$ (function over the tangent bundle) with respect to
supermap ${\rm{d}}\Phi$ (differential of the supermap $\Phi$).\\
The "inner product" of (pseudo)differential forms over $\Rm$ is a necessary
tool for building a "reasonable" physical theory on such supermanifold. If $f\,,h$
are homogenous elements with respect to Grassmann variables from the superspace
$[\C(\Pi T\Rm)]^{(k)}$, then their "inner product" is defined by
\be\label{"inner product"}
\langle f,\, h\rangle_g:=\mbox{\textbf{I}}[f\wedge *_g h]\ .
\ee
It is clear that the "inner product" $\langle\,.\,,\,.\,\rangle_g$ is $\mbb{R}$-linear and, moreover,
\begin{eqnarray}
\langle f,\,h\rangle_g & = & (-1)^k \langle h,\, f\rangle_g (-1)^{\tilde{f}\tilde{h}}\ ,\nonumber\\
0 & = & \langle \mbox{Even in\ }\xi,\, \mbox{Odd in\ }\xi\rangle_g \ ,\\
\langle Q(f),\, h\rangle_g & = & (-1)^{\tilde{f}}\langle f,\, \delta_g(h)\rangle_g\ .\nonumber
\end{eqnarray}
Let us note that the {\textbf{\emph{Laplace-DeRham}}} operator defined on (pseudo)differential forms
by
\be\label{Laplace-DeRham}
\Delta_g:=-(Q\delta_g+\delta_gQ)
\ee
is self-adjoint with respect to the "inner product" (\ref{"inner product"}). We use the
name "inner product" in quotation marks \emph{ex industria}, to emphasize the peculiar fact that
$\langle\,.\,,\,.\,\rangle_g$, given by (\ref{"inner product"}), is not always non-degenerate
(e.g. for the closed even 2-forms over supermanifold $\mbb{R}^{0|2}$ with metric
$g=\ep_{\a\b}\wedge \sigma^\a\wedge\sigma^\b$ it holds $\langle\,.\,,\,.\,\rangle_g\equiv 0)$.

\section{Grassmann Electrodynamics}

The symplectic mechanics and classical electrodynamics are undoubtedly nice and simple
applications of differential geometry in classical physics. The Grassmann electrodynamics is
a natural "grassmannisation" of the well known version of ordinary even electrodynamics to the
odd one. The inspiration for such a bit extravagant theory is provided by Cartan calculus
on a pure odd Riemannian supermanifold $\Rm$ and by the geometrical formulation of classical
electrodynamics (see e.g., interesting monographs \cite{Burke},\cite{Frankel}, or the
textbook \cite{Fecko}).\\
The electromagnetic field on the supermanifold $\Rm$ ("Grassmann space-time") is described
by the {\emph{potential 1-form}}
\be
A=A(\xi,y)=\mcal{A}_{\a}(\xi)y^\a\equiv\mcal{A}_{\a}(\xi)\,{\rm d}\xi^\a\ , \mbox{\ \ where\ \ } \mcal{A}_\a(\xi)\in[\bigwedge\mbb{R}^{2n}]_{[1]}\ .
\ee
Since the space of all odd potential $1$-forms on the supermanifold $\Rm$ is finite dimensional,
$\dim([\Pol(\Pi T\Rm)]^{(1)}_{[1]})=n2^{2n}$, the "Grassmann electromagnetism" possesses
only finite number of degrees of freedom.\\
The closed {\emph{electrodynamics $2$-form}}
\be
F:=Q(A)\in[\Pol(\Pi T\Rm)]^{(2)}_{[0]}\ ,
\ee
is invariant under the \emph{gauge transformations}}
\be\label{gauge transformation}
A\mapsto A^\prime=A+Q(f)\ , \mbox{\ \ where\ \ } f\in[\bigwedge\mbb{R}^{2n}]_{[0]}\ .
\ee
Externals sources in such theory are described by the {\emph{current $1$-from}}
\be\label{current}
J=J(\xi,y)=\mcal{J}_\a(\xi)y^\a\equiv\mcal{J}_\a(\xi)\,{\rm d}\xi^\a\ ,\mbox{\ \ where\ \ } \mcal{J}_\a(\xi)\in[\bigwedge\mbb{R}^{2n}]_{[1]}\ .
\ee
The dynamics of the \emph{Grassmann electromagnetic field} is governed by the action
\be\label{action}
S[A,\,J,\,g]:=-\frac{1}{2}\langle Q(A)\,,\,Q(A)\rangle_g+\langle J\,,\,A \rangle_g\ ,
\ee
which is a functional (strictly speaking only a function, because the Grassmann electrodynamics has
finite number of degrees of freedom) of the potential and current $1$-forms $A$ and $J$ and also of
the metric tensor $g$. Variation of the action (\ref{action}) with respect to the potential $A$
(without any boundary condition) leads (in the case when $\langle\,.\,,\,.\,\rangle_g$ is regular
"inner product" over 1-forms) to the set of algebraic {\emph{Grassmann-Maxwell equations}},
which can be written in a compact form as
\be\label{Grassmann-Maxwell}
\delta_g F=-J\ .
\ee
The Grassmann-Maxwell equations are internally consistent only if the current $1$-form
$J$ is co-closed, i.e. if
\be\label{continuity equation I}
\delta_g J=0\ \ \Leftrightarrow \ \ Q(*_g\, J)=0\ .
\ee
This is the {\emph{continuity equation}} for the external currents.\\
To bring this extravagancy to the top of its bent, we split coordinates covering "Grassmann space-time"
$\Rm=\mbb{R}^{0|(2n-1)}\oplus\Rz$ to the "space" coordinates
($\xi^1,\dots ,\xi^{2n-1}$) and the "time" coordinate ($\xi^{2n}=\tau$).
We restrict the closed integral form $*_g J\in *(\Pol(\Pi T\Rm))$ to the "space" hyper-surface
$\mbb{R}^{0|(2n-1)}$. The \emph{Grassmann charge} related to external sources is then defined as
an integral of the restricted form $(*_g J)\mid_{{\xi^{2n}=\tau \atop y^{2n}=0}}$
over "Grassmann space" subsupermanifold $\mbb{R}^{0|(2n-1)}$, namely
\be
\mbox{\textbf{Q}}_{Gr}(\tau):=\int\limits_{\bigwedge\mbb{R}^{2n-1}}\overline{{\rm d}\xi} \int\limits_{\mbb{R}^{2n-1}}\overline{{\rm d}y}\,\biggl.(*_g J)(\xi,y)\biggr|_{{\xi^{2n}=\tau \atop y^{2n}=0}}\equiv\mbox{\textbf{I}}_{Sp}\biggl[(*_g J)\mid_{{\xi^{2n}=\tau \atop y^{2n}=0}}\biggr]\ .
\ee
The odd vector field $\mcal{T}=\pd_\tau$ defines a corresponding "Grassmann time" evolution. It is clear
from equations (\ref{Lie derivation}), (\ref{Cartan magic formula}) and (\ref{continuity equation I}) that
\bdis
\mcal{T}(*_g J)=\mcal{T}^\ua(*_g J)=Q[\mcal{T}_\ua(*_g J)]+\mcal{T}_\ua[Q(*_g J)]=Q[\mcal{T}_\ua(*_g J)]\ .
\edis
Last equation and the odd \emph{Stokes theorem} (\ref{per-partes & Stokes theorem}) imply that
\be
0=-\mbox{\textbf{I}}_{Sp}\biggl[Q\biggl(\mcal{T}_\ua(*_g J)\mid_{{\xi^{2n}=\tau \atop y^{2n}=0}}\biggr)\biggr]=\cdots=\mcal{T}\biggl(\mbox{\textbf{I}}_{Sp}\biggl[(*_g J)\mid_{{\xi^{2n}=\tau \atop y^{2n}=0}}\biggl]\biggr)=\frac{\pd\, \mbox{\textbf{Q}}_{Gr}}{\pd\, \tau}\ .
\ee
Let us remark that similarly as in the ordinary electromagnetism, it is
possible to fix gauge freedom (of electromagnetic 2-form $F$) by imposing the condition
\be\label{Lorentz gauge condition}
\delta_g A=0\ .
\ee
This {\emph{Lorentz gauge condition}} does not define the function $f$ in the gauge transformation
(\ref{gauge transformation}) uniquely. Let us note that this uncertainty in the choice of the gauge fixing
function $f$ is similar as in the ordinary electrodynamics: to any solution of the Poisson equation
$\Delta_g f=\delta_g A$, which determines $f$, it is always possible to add a solution of the homogeneous
(Laplace) equation.\\
To eliminate this ambiguity in the ordinary electromagnetism ($g_{\a\b}=\eta_{\a\b}$, i.e. in $(3+1)$-Minkowski space-time),
we need to fix the boundary conditions for the gauge fixing function $f$ (usually $f(x)\rightarrow 0$
for $x\rightarrow\infty$), because then the Poisson equation $\Delta_\eta f=\delta_\eta A$ has unique
solution\footnote{Let us note that
this statement is not true for the general Riemannian manifold $(M,\,g)$, since the solution of the Poisson equation
in particular physical situations constitutes a serious problem of mathematical physics.}.\\
Unfortunately, this does not work in the case of the Grassmann electrodynamics, where the boundary
condition could be most naturally established by requirement that $0=\left.f(\xi)\right|_{\xi=0}$
(the absolute term in $f$ does not contribute to gauge transformation (\ref{gauge transformation})). The problem
consists in an unpleasant fact that the \emph{Laplace-DeRham operator} restricted to the functions\footnote{sometimes
called as the \emph{Laplace-Beltrami operator}}
\be
\Delta_g f=\delta_g\,Q (f)=\frac{\pd\ }{\pd\,\xi^\a}\,g^{\a\b}\frac{\pd f}{\pd\,\xi^\b}-\frac{\pd\ln{\sqrt{|g|}}}{\pd\,\xi^\a}\,g^{\a\b}\frac{\pd f}{\pd\,\xi^\b}
\ee
is not invertible. This peculiar fact has a serious consequence: it is not possible
to quantize the Grassmann-Maxwell theory by using Faddeev-Popov method in a straightforward way,
because after applying gauge
fixing condition (\ref{Lorentz gauge condition}) there still remains relative "rich" gauge freedom, which is
represented by the subgroup $\mbox{\it{Ker\,}}\Delta_g$. Similar problems arise also with another
gauge fixing conditions, e.g. $\xi^\mu\mcal{A}_\mu(\xi)=0$ (\emph{radial gauge}) or $\mcal{A}_{2n}(\xi)=0$
(\emph{Hamiltonian gauge}) and therefore the problem of the quantization will not be discussed in this paper.
The reason is down-to-earth, needs an individual approach from the case to case.\\
The Grassmann-Maxwell equations (\ref{Grassmann-Maxwell}) could be rewritten in Lorentz gauge
(\ref{Lorentz gauge condition}) in the form
\be
\Delta_g A=J\ .
\ee
Now, let us consider the case of a free electromagnetic field ($J\equiv 0$). The variation of the action
(\ref{action}) with respect to the metric $g$
\be\label{variation under metric}
S[A,\,g+\Delta t h]=S[A,\,g]-\Delta t\int\limits_{\bigwedge\mbb{R}^{2n}}\frac{\overline{{\rm d}\xi}}{\sqrt{|g|}}\,\frac{1}{2}\,h_{\a\b}(\xi)\,T^{\a\b}(\xi)+o(\Delta t)\ ,
\ee
where $\Delta th$ is the perturbation of metric ($\Delta t$ is an infinitesimal real parameter)
and $\overline{{\rm d}\xi}\sqrt{|g^{-1}|}$ is invariant Berezin integral measure, defines
the covariant tensor field (a quadratic function in the fiber variables of the tangent bundle):
\be
T=T(\xi,\,\sigma)=T_{\a\b}(\xi)\wedge\sigma^\a\wedge\sigma^\b\,,\ \ \ \mbox{where}\ \ \ T_{\a\b}:=g_{\a\mu}g_{\b\nu}T^{\mu\nu}\ .
\ee
It is the {\emph{energy-momentum tensor}} of the free Grassmann electromagnetic field. It is clear from
definition equation (\ref{variation under metric}) that in the energy-momentum tensor
$T$ only the skew-symmetric part is uniquely fixed (contrary to standard electrodynamics, where
the same is true of the symmetric part). Namely,
{\small
{\bdis\label{energy-momentum tensor}
T_{\a\b}=-\frac{1}{2}\int\limits_{\mbb{R}^{2n}}\int\limits_{\mbb{R}^{2n}}\frac{{\overline{{\rm d}y}}{\overline{{\rm d}z}}}{\ (2\pi)^n}|g|\wedge\{g_{\a\b}+\imath(z_\a y_\b-z_\b y_\a)\}\wedge\{F(\xi,y)\wedge F(\xi,z)\}\wedge e^{-\imath z^\mu g_{\mu\nu} y^\nu}\ ,
\edis}}%
where $z_\a=g_{\a\mu}z^\mu$ and $y_\b=g_{\b\nu}y^\nu$.\\

\noindent The energy-momentum tensor is a fundamental physical object that enables us
to construct preserving quantities, which correspond to the space-time symmetries associated
with (conformal) Killing's vector fields. Our next steps will be more general and final
results could be automatically applied to Grassmann electrodynamics discussed above.\\
Consider now an arbitrary field $\psi$ (function, 1-form, tensor field,...)
over the supermanifold $\Rm$ with fixed metric $g$. Let us assume that the dynamics of the field $\psi$ is governed
by the action
\be\label{generalized action}
S[\psi,\,g]=\int\limits_{\bigwedge\mbb{R}^{2n}}\frac{\overline{{\rm d}\xi}}{\sqrt{|g|}}\ \biggl[L(\psi,\,g)\biggr](\xi)\ ,
\ee
which is natural with respect to diffeomorphisms, i.e.
\be
\Phi^*\biggl(L(\psi,\,g)\biggr)=L\biggl(\Phi^*(\psi),\,\Phi^*(g)\biggr)\ ,
\ee
for any $\Phi\in\Diff(\Rm)$. If $V\in\X(\Rm)$ is an arbitrary even vector field, $V$ defines the
infinitesimal flow $\Phi_{^{\Delta t}}(V)$ on the Riemannian supermanifold $\Rm$. The associated
horizontally lifted even vector fields $V^\ua,\ V^{\ua_{hor}},...$ generate the corresponding flows
$\Phi_{^{\Delta t}}(V^{\ua_{...}})$ over supermanifolds $\Pi T\Rm,\ T\Rm,...$ and consequently
the \emph{pull-backs} $\Phi^*_{^{\Delta t}}(V^{\ua_{...}})$ over forms, symmetric covariant tensors fields and
so on. Since the berezinians of all such lifted flow-transformations are identically equal to unity and because
the action (\ref{generalized action}) is natural with respect to supergroup $\Diff(\Rm)$, we could write
\begin{eqnarray*}
S[\psi,\,g] & = & S[\Phi^*_{^{\Delta t}}(V^{\ua_{...}})(\psi),\,\Phi^*_{^{\Delta t}}(V^{\ua_{hor}})(g)] \\
{}          & = & S[\psi+\Delta t V^{\ua_{...}}(\psi)+o(\Delta t),\,g+\Delta t V^{\ua_{hor}}(g)+o(\Delta t)]\ .
\end{eqnarray*}
If the field $\psi$ is the solution of the dynamical equations, i.e. $\frac{\delta S}{\delta \psi}=0$,
previous equation implies that
\be\label{equation}
0=\int\limits_{\bigwedge\mbb{R}^{2n}}\frac{\overline{{\rm d}\xi}}{\sqrt{|g|}}\ \frac{1}{2}(V^{\ua_{hor}}g)_{\a\b}\,T^{\a\b}\ ,\ \ \ \mbox{where}\ \ \ T^{\a\b}(\xi)=-2\sqrt{|g|}\,\frac{\delta S[\psi,\,g]}{\delta g_{\a\b}(\xi)}\ .
\ee
Explicit coordinate expression (\ref{horizontal lift}) for the horizontal lift $V^{\ua_{hor}}$ ($V$ is even) allows us
to write
\be
\frac{1}{2}(V^{\ua_{hor}}g)_{\a\b}\,T^{\a\b}=V^\mu\Biggl[\frac{1}{2}\,\frac{(g_{\a\b})_{,\,\mu}T^{\a\b}}{\sqrt{|g|}}-\Biggl(\frac{g_{\mu\a}T^{\a\b}}{\sqrt{|g|}}\Biggr)_{,\,\b}\Biggr]-\Biggl(\frac{V^\mu g_{\mu\a} T^{\a\b}}{\sqrt{|g|}}\Biggr)_{,\,\b}\ .
\ee
Since the equation (\ref{equation}) is valid for all even vector fields $V\in\X(\Rm)$ (components
$V^\mu(\xi)=V(\xi^\mu)$ are odd functions and $\int\overline{{\rm d}\xi}\,(\ )_{,\,\b}\equiv 0$), we conclude that
\be\label{continuity equation}
\mbox{$\psi$ satisfies\ } \frac{\delta S}{\delta \psi}=0\ \ \Rightarrow\ \ \frac{1}{2}\,\frac{(g_{\a\b})_{,\,\mu}T^{\a\b}}{\sqrt{|g|}}-\biggl(\frac{g_{\mu\a}T^{\a\b}}{\sqrt{|g|}}\biggr)_{,\,\b}\equiv 0\ .
\ee
This equation is a compact formulation of \emph{conservation laws}, because from ordinary point of view,
equation (\ref{continuity equation}) is equivalent to the statement that $(T^{\mu\b})_{;\,\b}=0$.
Now, it becomes clear that for an arbitrary Killing's vector field $V\in\X(\Rm)$ (the generator of an isometry of $\Rm$)
the 1-form
\be
T_V=T_V(\xi,\,y):=\frac{1}{2}\Biggl. V_{\ua_{ver}}\biggl(T_{\a\b}(\xi)\wedge\sigma^\a\wedge\sigma^b\biggr)\Biggr|_{\sigma\,\mapsto\, y}\
\ee
is co-closed, i.e. $Q(*_g T_V)=0$ (occasionally we are using metric to raise and lower indices,
in order to make our final formulae more compact). If moreover $T^\mu_{\ \mu}=0$, the same
is true for the arbitrary conformal Killing's vector field and therefore the quantity
\be
\mbox{\textbf{Q}}_{V}(\tau):=\mbox{\textbf{I}}_{Sp}\biggl[(*_g T_V)\mid_{{\xi^{2n}=\tau \atop y^{2n}=0}}\biggr]
\ee
is an integral of motion (i.e. $\pd_{\tau}[\mbox{\textbf{Q}}_{V}(\tau)]=0$) and $T_V$ is corresponding
\emph{Noether's current} associated to the (conformal) Killing's vector field $V\in\X(\Rm)$.
The proof of this assertion is the same as for the \emph{Grassmann charge} $\mbox{\textbf{Q}}_{Gr}$ discussed above.

\section{Connection and Grassmann General Relativity}

In this section we shall try to uncover the Grassmann analogy of the Einstein's theory of General Relativity. We shall
derive equations, which will govern the dynamics of the metric on the pure odd Riemannian supermanifold
$(\Rm\,,\,g)$, but before this we are forced to establish a necessary supermathematical tool, namely the linear connection.
The inspiration for such "monstrosity" comes from the standard General Relativity (see e.g., \cite{Weinberg},
\cite{Landau-Lifshitz II}) and from the supergeometry discussed above.\\
We start with the definition of the \emph{\textbf{linear connection}} (the quintessence of general relativity)
on functions and vector fields for the general (not necessary Riemannian) pure odd supermanifolds. Connection
(covariant derivative $\nabla$ with respect to vector field $V\in\X(\R)$) over functions $\mcal{F}(\R)$ is trivial:
we require (analogically as in the ordinary case) that $\nabla_V f:=V f$. Connection over vector fields
is the operation
\be\label{connection}
\nabla:\ \biggl(\X(\R)\biggr)^2\rightarrow\X(\R)\ ,\ \ \ \ \ \ \ \ \ \ \ \  (V,\,W)\mapsto\nabla_V W\ ,
\ee
which in addition for arbitrary vector fields $V,\,W,\,U\in\X(\R)$ and function $f\in\mcal{F}(\R)$ (all this objects
are considered here and below as homogenous elements with respect to parity) satisfies:
\begin{eqnarray}\label{definition of connection}
\widetilde{\nabla_V W} & = & \tilde{V}+\tilde{W}\ ,\nonumber\\
\nabla_{V+fU}W & = & \nabla_{V}W+f\nabla_{U}W\ ,\\
\nabla_V (W+fU) & = & \nabla_V W +(V f)U+(-1)^{\tilde{V}\tilde{f}}f\,\nabla_V U \ .\nonumber
\end{eqnarray}
It is clear that the full information about connection $\nabla$ on supermanifold $\R$ is uniquely encoded into
$n^3$ Christoffel symbols (odd functions) $\Gamma^{\g}_{\a\b}:=(\nabla_{\pd_{\xi^\b}}\pd_{\xi^{\a}})(\xi^\g)$,
i.e. $\Gamma^{\g}_{\a\b}(\xi)$ are components of the vector field $\nabla_{\pd_{\xi^\b}}\pd_{\xi^{\a}}$ with respect
to coordinates frame $(\xi^1,\dots ,\xi^n)$. If $\xi^\a\mapsto\Xi^\a(\xi)$, then $\Gamma$'s transform as
\bdis
{\Gamma^\prime}^{\g}_{\a\b}=\frac{\pd\, \Xi^\g}{\pd\, \xi^\l}\,\frac{\pd\, \xi^\m}{\pd\, \Xi^\a}\,\frac{\pd\, \xi^\n}{\pd\, \Xi^\b}\,\Gamma^\l_{\m\n}+\frac{\pd\, \xi^\m}{\pd\, \Xi^\a}\,\frac{\pd\, \xi^\n}{\pd\, \Xi^\b}\,\frac{\pd\ }{\pd\, \xi^\mu}\,\biggl(\frac{\pd\, \Xi^\g}{\pd\, \xi^\n}\biggr)\ .
\edis
In the same spirit as in the case of ordinary connection, the \emph{\textbf{torsion}} associated to connection
$\nabla$ is the (super)tensorial (i.e. graded skew-symmetric and graded $f$-multilinear) operation $\T_n$
defined by
\begin{eqnarray}\label{torsion}
\T_n\ :\biggl(\X(\R)\biggr)^2 & \rightarrow & \X(\R)\ , \\
(V,W)\hspace{0.5cm} & \mapsto & \T_n(V,W):=\nabla_VW-(-1)^{\tilde{V}\tilde{W}}\nabla_WV-[V,W]\ .\nonumber
\end{eqnarray}
Let us note that $\T_n(\pd_{\xi^\a},\pd_{\xi^\b})=(\Gamma^{\g}_{\a\b}+\Gamma^{\g}_{\b\a})\pd_{\xi^\g}
=\T_n(\pd_{\xi^\b},\pd_{\xi^\a})$ and therefore the torsion vanishes identically only in the case, when the
connection $\nabla$ on the supermanifold $\R$ is anti-symmetric, i.e. $\Gamma^{\g}_{\a\b}=-\Gamma^{\g}_{\b\a}$.\\
The \emph{\textbf{curvature tensor}} $\U_n$ assigned to the connection $\nabla$ is naturally defined as the supermap
\begin{eqnarray}\label{curvature}
\U_n:\ \biggl(\X(\R)\biggr)^3 & \rightarrow & \X(\R)\ ,\\
(V,W,U)\hspace{0.27cm} & \mapsto & \U_n(V,W,U):=[\nabla_V,\nabla_W]U-\nabla_{[V,W]}U\nonumber\ ,
\end{eqnarray}
where $[\nabla_V,\nabla_W]=\nabla_V\nabla_W-(-1)^{\tilde{V}\tilde{W}}\nabla_W\nabla_V $.
Obvious coordinate expression for the curvature tensor $\U_n$ is
\be\label{components of curvature tensor}
\U_n(\pd_{\xi^\m},\pd_{\xi^\n},\pd_{\xi^\b})=\mcal{R}^{\a}_{\b\m\n}\pd_{\xi^\a}=\biggl\{(\Gamma^{\a}_{\b\m})_{,\,\n}+(\Gamma^{\a}_{\b\n})_{,\,\m}-\Gamma^{\l}_{\b\m}\Gamma^{\a}_{\l\n}-\Gamma^{\l}_{\b\n}\Gamma^{\a}_{\l\m}\biggr\}\pd_{\xi^\a}\ .
\ee
Apart from the basic properties of the torsion $\T_n$ and curvature $\U_n$, which have arisen as straightforward
consequences of their definitions, one also obtains two useful and powerful identities, namely
the \emph{Ricci identity}
\be\label{Ricci identity}
\sum\limits_{cycl.\atop permut.}(-1)^{\tilde{V}\tilde{U}}\Biggl\{\U_n(V,W,U)-\nabla_V \biggl(\T_n(W,U)\biggr)-\T_n(V,[W,U])\Biggr\}=0
\ee
and the \emph{Bianchi identity}
\be\label{Bianchi identity}
\sum\limits_{cycl. \atop permut.}(-1)^{\tilde{V}\tilde{U}}\Biggl\{\biggl[\nabla_V,\U_n(W,U,.)\biggr]+\U_n(V,[W,U],.)\Biggr\}=0\ .
\ee
The proof of these assertions is not difficult, it is mainly a technical matter, which needs patience, large
sheet of paper and knowledge about the (super)Jacobi identity stating that
$(-1)^{\tilde{V}\tilde{U}}[V,[W,U]]+\mbox{cycl.\ permut.}=0$ for any $V,\,W,\,U\in\X(\R)$.\\
On the odd Riemannian supermanifold $(\Rm,\,g)$ it is possible to require, analogically as in the ordinary case,
the compatibility between connection $\nabla$ and metric $g$, which means that for arbitrary vector fields
$V,\,W,\,U$ on $\Rm$ it should holds
\be\label{compatibility with metric I}
\nabla_V(W\,,\,U)_g=(\nabla_V W\,,\,U)_g+(-1)^{\tilde{V}(\tilde{W}+1)}(W\,,\,\nabla_V U)_g\ .
\ee
As was mentioned above, the connection $\nabla$ can be uniquely defined by $(2n)^3$-odd functions ($\Gamma$),
compatibility condition (\ref{compatibility with metric I}) gives $2n^2(2n-1)$-secondary constrains on the
Christoffel symbols
\bdis
(g_{\a\b})_{,\, \g}=g_{\a\m}\Gamma^{\m}_{\b\g}-g_{\b\m}\Gamma^{\m}_{\a\g}\ .
\edis
But there still remains $2n^2(2n+1)$-degrees of freedom for the Christoffel $\Gamma$'s, the residual ambiguity
could be eliminated by requirement that connection $\nabla$ is without torsion ($\T_n=0$) or equivalently that
$\nabla$ is anti-symmetric.
The final expression for the Christoffel symbols related to the \emph{metric and anti-symmetric connection} $\nabla$
is very simple, and similar to the ordinary one. Namely,
\be
\Gamma_{\a\b}^\g=\frac{1}{2}\,g^{\g\m}\,\biggl[(g_{\m\a})_{,\, \b}-(g_{\m\b})_{,\, \a}-(g_{\a\b})_{,\, \m}\biggr]\equiv g^{\g\m}\,\Gamma_{\m\a\b}\ .
\ee
The pure metricity of the connection $\nabla$ (in principle $\T_n\neq 0$) implies that
\be\label{symmetry I}
\biggl(\U_n(V,W)X\,,\,Y\biggr)_g=(-1)^{(\tilde{X}+1)(\tilde{Y}+1)}\biggl(\U_n(V,W)Y\,,\,X\biggr)_g\ ,
\ee
the last equation together with anti-symmetry of $\nabla$ reveal another interesting symmetry of the
supermap $\U_n$:
\be\label{Symmetry II}
\biggl(\U_n(V,W)X\,,\,Y\biggr)_g(-1)^{\tilde{Y}}=(-1)^{(\tilde{V}+\tilde{W})(\tilde{X}+\tilde{Y})}\biggl(\U_n(X,Y)V\,,\,W\biggr)_g(-1)^{\tilde{W}}\ .
\ee
The proof of these identities is very simple and it is primarily based on the fact that for any function
$f\in\mcal{F}(\Rm)$ it holds $[\nabla_V,\nabla_W]f-\nabla_{[V,W]}f=0$. The choice $f=(X\,,\,Y)_g$
gives (\ref{symmetry I}). The statement (\ref{Symmetry II}) is just a consequence of previous one together with Ricci
identity without torsion terms.\\
Henceforth, we work only with the \emph{Grassmann-Levi-Civita} $\nabla$, which is torsionless and compatible with
metric $g$. In such case, it is evident that
the curvature tensor $\U_n$ over the Riemannian supermanifold $(\Rm,\,g)$ has only
$\frac{n^2}{3}[4n^2-1]+n[4n^2+1]$ independent components $\mcal{R}^{\a}_{\b\m\n}$, which is indeed
more than in ordinary case for the same dimension.\\

\noindent Now we are able to describe the gravity on $\Rm$. The free dynamics of the metric $g$ is governed
by the \emph{Grassmann-Hilbert action}
\be\label{Grassmann-Hilbert action}
S_{G-H}[g]:=\int\limits_{\bigwedge\mbb{R}^{2n}}\frac{\overline{{\rm d}\xi}}{\sqrt{|g|}}\ \mcal{R}\ ,
\ee
where the \emph{Ricci scalar} $\mcal{R}$ is defined by
\be\label{Ricci scalar}
\mcal{R}:=g^{\b\a}\mcal{R}^{\m}_{\a\m\b}=g^{\b\a}\biggl\{(\Gamma^{\m}_{\a\b})_{,\,\m}+(\Gamma^{\m}_{\a\m})_{,\,\b}-\Gamma^{\l}_{\a\b}\Gamma^{\m}_{\l\m}-\Gamma^{\l}_{\a\m}\Gamma^{\m}_{\l\b}\biggr\}\ .
\ee
The first two terms in the action (\ref{Grassmann-Hilbert action}) can be simplified by \emph{per-partes}
integration ($\int\overline{{\rm d}\xi}\,(\ )_{,\,\m}\equiv 0$).
After such treatment we obtain more convenient expression, namely
\be\label{simple G-H action}
S_{G-H}[g]=\int\limits_{\bigwedge\mbb{R}^{2n}}\frac{\overline{{\rm d}\xi}}{\sqrt{|g|}}\ \mcal{G}\ ,\mbox{\ \ \ where\ \ \ } \mcal{G}=g^{\b\a}\,\biggl\{\Gamma_{\a\b}^\m\Gamma_{\m\l}^\l+\Gamma^\l_{\a\m}\Gamma^{\m}_{\l\b}\biggr\}\ .
\ee
The variation of (\ref{simple G-H action}) with respect to the metric $g$ (without any boundary conditions, forasmuch as
the surface terms do not contribute)
leads to the system of algebraic equations
\be
\frac{\pd (\sqrt{|g^{-1}|}\,\mcal{G})}{\pd g_{\a\b}}-\frac{\pd\ }{\pd \xi^\m}\Biggl[\frac{\pd(\sqrt{|g^{-1}|}\,\mcal{G})}{\pd (g_{\a\b})_{,\,\m}}\Biggr]=0\ ,
\ee
which dictates the evolution of the free gravitational field. It is always possible to add to the \emph{Grassmann-Hilbert}
action a mater term of the form (\ref{generalized action}), and then the variation of
$S_{G-H}[g]+S_{mater}[\psi,\,g]$ with respect to $g$, gives the \emph{Grassmann-Einstein equations}
\be
\frac{\pd (\sqrt{|g^{-1}|}\,\mcal{G})}{\pd g_{\a\b}}-\frac{\pd\ }{\pd \xi^\m}\Biggl[\frac{\pd(\sqrt{|g^{-1}|}\,\mcal{G})}{\pd (g_{\a\b})_{,\,\m}}\Biggr]=\frac{1}{2}\sqrt{|g^{-1}|}\,T^{\a\b}\ ,
\ee
where $T^{\a\b}$ is the energy-momentum tensor (\ref{equation}) related with the mater fields.

\begin{center} {\textbf{Acknowledgement}} \end{center}

\noindent I would like to express my deep gratitude to my teacher P. \v Severa,
who was an excellent guide of mine to supergeometry. I am also very grateful to P.
Pre\v snajder and V. Balek who gave me much valuable advice and who carefully read
a draft of this paper and many thanks go to my girlfriend Tulka ($\heartsuit$) for her
encouragement.

\end{document}